\documentclass[12pt]{amsart}
\usepackage{a4wide}
\usepackage{amssymb}
\usepackage{mathrsfs}
\theoremstyle{plain}
\newtheorem{theorem}{Theorem}
\newtheorem{lemma}{Lemma}[section]
\newtheorem{proposition}{Proposition}[section]

\newcommand{\NN}{\mathbf{N}}
\newcommand{\N}{\mathbf{N}} \newcommand{\M}{\mathbf{M}} 
\newcommand{\bR}{\mathbf{R}}

\newcommand{\e}{\epsilon}
\newcommand{\ZZ}{\mathbf{Z}} \newcommand{\BB}{\mathbf{B}}
\newcommand{\uZZ}{\underline{\ZZ}}
\newcommand{\ff}{\mathbf{f}} 
\newcommand{\bg}{\mathbf{g}}\newcommand{\bh}{\mathbf{h}}
\newcommand{\bd}{\mathbf{d}}

\newcommand{\hh}{\mathbf{h}}
\newcommand{\C}{\mathbf{C}}
\renewcommand{\P}{\mathbf{P}}

\newcommand{\X}{\mathbf{X}} 
\newcommand{\bX}{\underline{\X}}\newcommand{\bx}{\mathbf{x}}
\newcommand{\bux}{\underline{\mathbf{x}}}
\newcommand{\by}{\mathbf{y}} \newcommand{\buy}{\underline{\mathbf{y}}}
\newcommand{\A}{\mathbf{A}}\newcommand{\bS}{\mathbf{S}}
\newcommand{\bA}{\underline{\A}}

\newcommand{\cL}{\mathcal{L}} 
\newcommand{\cA}{\mathcal{A}} \newcommand{\cB}{\mathcal{B}}
\newcommand{\cC}{\mathcal{C}}
 
\newcommand{\cP}{\mathcal{P}} \newcommand{\cM}{\mathcal{M}}
\newcommand{\ubox}{\underline{\square}_1}
\newcommand{\pp}{\mathbf{p}}
\newcounter{exa}
\newcommand{\hlimo} {\hat{lim}_\omega}
\newcommand{\limo} {\lim_\omega}
\begin{document}

\title[]{Samplings and observables. Invariants of metric measure
spaces.}
\author[G. Elek]{G\' abor Elek}
\thanks{Work supported in part 
by a Marie Curie grant, TAMOP 4.2.1/B-09/1/KMR-2010-003 and MTA Renyi "Lendulet" Groups and
Graphs Research Group.}
\begin{abstract} In the paper we are dealing with metric measure spaces
of diameter at most one and of total measure one. Gromov introduced the 
sampling compactification of the set of these spaces. He asked whether
the metric measure space invariants extend to the compactification.
Using ideas of the newly developed theory of graph limits we identify
the elements of the compactification with certain geometric objects  and
show how to extend various invariants to this space. We will introduce
the notion of ultralimit of metric measure spaces, that will be the main
technical tool of our paper.

\end{abstract} 
\maketitle
\tableofcontents
\section{Introduction and preliminaries}
\label{intro}
\subsection{Metric measure spaces}
Let $\chi$ denote the set (up to isomorphisms) of 
Polish spaces $X$ of diameter at most one equipped
with a Borel probability measure $\mu_X$ of full support. In the course
of the paper, we refer to these objects as mm-spaces.
We say that $(X,\mu_X)$ and
$(Y,\mu_Y)$ are equivalent if there exists an isometry 
$\Phi:X\to Y$ such that
$\Phi_*(\mu_X)=\mu_Y$. The space $\chi$ has a Polish space structure.
Note that we will often use the notation $X$ instead of $(X,\mu_X)$. 
The distance of $X,Y\in\chi$ is defined the following way \cite{Gro}.
Consider the set $MPM(X)$
of measure preserving maps $\Psi ([0,1],\lambda)\to (X,\mu_X)$, where
$\lambda$ is the Lebesgue measure. Let $\Psi^{-1}(d_X)$ be the pull-back
of the distance function of $X^2$. Then
$$\ubox(X,Y)=\inf_{\Psi_1\in MPM(X),\Psi_2\in MPM(Y)}
\square_1(\Psi_1^{-1}(d_X),\Psi_2^{-1}(d_Y))\,,
$$ where the $\square_1$-distance of the measurable functions 
$f,g:[0,1]^2\to [0,1]$ is defined as the supremal $\e$ such that
$f$ and $g$ are $\e$-close outside a subset $X_\e\subset [0,1]$ of
measure at most $\e$ that is,
$$|f(x,y)-g(x,y)|\leq\e$$
for $(x,y)\in ([0,1]\backslash X_\e)\times ([0,1]\backslash X_\e)$.
\subsection{Samplings}
Gromov introduced an other notion of convergence in $\chi$, convergence
in samplings \cite{Gro} (see also \cite{V}). 
Let $M_\infty$ be the convex compact metric
space of infinite matrices $\{d_{i,j}\}^\infty_{i,j\geq 1}$
with $0\leq d_{i,j}\leq 1$, $d_{i,j}=d_{j,i}, d_{i,i}=0\,.$
We have a continuous
map $\rho:X^{\NN}\to M_\infty$ that assigns $\{d_{i,j}=d_X(x_i,x_j)\}$ to
the sequence $(x_1,x_2,\dots)\in X^{\NN}$. We denote by $\mu^X_\infty$ the
push-forward of the measure $\mu_X\times\mu_X\dots$.
According to the Reconstruction Theorem ($3^{\frac{1}{2}} 5.$ \cite{Gro}),
the map $\tau:\chi\to \cP(M_\infty)$ is a continuous injective map
(where $\cP(M_\infty)$
is the space of probability measures on $M_\infty$).
Thus the closure of $\tau(\chi)$ gives us a compactification 
$\overline{\chi}$ of 
$\chi$.
We say that $\{X_n\}^\infty_{n=1}$ is convergent in sampling if
$\{\tau(X_i)\}^\infty_{i=1}$ is a weakly convergent sequence.
We use the notation $X_n\stackrel{s}{\to} X$, if
$\{\tau(X_i)\}^\infty\to\tau(X)$ weakly.
One of the main motivation for writing this paper was a remark of
Gromov in Section 3$\frac{1}{2}$ of his book \cite{Gro}. He asked whether
the invariants of mm-spaces can be extended to $\overline{\chi}$.
In our paper we argue that the answer is yes, the elements of $\overline{\chi}$
can be regarded as geometric objects.

\noindent
Let $G$ be a finite simple graph. We can associate an element $X_G\in\chi$
the following way. $X_G=V(G)$, $d_X(a,b)=1/2$ if $a$ and $b$ are connected
and $d_X(a,b)=1$ otherwise.
It is important
to note that $\{G_n\}^\infty_{n=1}$ is convergent as a dense graph sequence
(see e.g. \cite {BCLV})
if and only if $\{X_{G_n}\}^\infty_{n=1}$ is convergent in sampling. 
Our paper uses the ideas and methods of graph limit theory in a 
substantial way.

\subsection{Quantum metric spaces}
Let $(X,\mu)$ be standard Borel space with a probability measure and
$d^*:X\times X\to\cP[0,1]$ be a probability measure space valued 
Borel-measurable function, such
that $d^*(x,x)=\delta_0$ and $d^*(x,y)=d^*(y,x)$. Then for a triple of
different points $a,b,c$ points we can pick independent random lengths
$l(a,b), l(b,c), l(a,c)$ using the probability measures $d^*(a,b),
d^*(b,c), d^*(a,c)$. If for $\mu$-almost triples the random lengths
satisfy the triangle inequality, we call $(X,\mu,d^*)$ a quantum metric
measure space, a qmm-space.
Then, for any qmm-space $X$ we can associate an $S_\infty$-invariant
measure $\tau(X)$ the following way.
First consider $\hat{M}_\infty$ the convex compact metric space of infinite
matrices with coefficients $d*_{i,j}$ in $\cP[0,1]$ such that 
$d^*_{i,j}=d^*_{j,i}, d^*_{i,i}=\delta_0\,.$
An element $\nu\in\hat{M}_\infty$ defines a probability measure on
$M_\infty$. By the definition of qmm-spaces, we have
a continuous map $X^\N\to \hat{M}_\infty$ and we consider the push-forward
of $\mu_X\times\mu_X\times\dots$. 
This is a probability measure on a convex compact
space (i.e the space of probability measures on $M_\infty$) so one can
consider its barycenter. This will be the associated probability measure
$\tau((X,\mu,d))$. We will prove that qmm's completely represent
$\overline{\chi}$.
\begin{theorem} \label{elsotetel}
If $\kappa\in\overline{\chi}$ then there exists a qmm-space $X$
such that $\tau(X)=\kappa$.
\end{theorem}
We can extend the metric $\ubox$ to qmm-spaces as well. 
If $f$ and $g$ are $\cP[0,1]$-valued weak-$*$-measurable functions on $[0,1]^2$
then their $\square_1$-distances can be defined again as
the supremal supremal $\e$ such that
$f$ and $g$ are $\e$-close in the $d_{ext}$-distance, outside a subset 
$X_\e\subset [0,1]$ of
measure at most $\e$.
Here $d_{ext}$ is the metric extension of the usual distance on $[0,1]$
onto $\cP[0,1]$. That is
$$d_{ext}(\mu,\nu)=\sup_{f\in Lip_1[0,1]} |\int_0^1 f\,d\mu-\int_0^1 f\, d\nu|\,$$
Thus
$$\ubox(X,Y)=\inf_{\Psi_1\in MPM(X),\Psi_2\in MPM(Y)} 
\square_1((\Psi_1)^{-1}(d_X), (\Psi_2)^{-1}(d_Y))\,.$$
Note that two qmm-spaces with zero $\ubox$-distance are not
necessarily isomorphic. However, we will prove the following
reconstruction theorem.
\begin{theorem} \label{masodiktetel}
If $X,Y$ are qmm-spaces then $\tau(X)=\tau(Y)$ if and only if $\ubox(X,Y)=0$.
In fact, if $\tau(X)=\tau(Y)$, then there exist 
$\Psi_1\in MPM(X)$ and $\Psi_2\in MPM(Y)$ such
that $(\Psi_1)^{-1}(d_X)=(\Psi_2)^{-1}(d_Y)$.
\end{theorem}
It is important to note that the idea of qmm-spaces is already
implicit in the work of Lov\'asz and Szegedy \cite{LSZDEC}. 
Also, Theorem \ref{masodiktetel} is an analogue of the
uniqueness theorem of \cite{BCL}, proved originally for graphons (which
can be regarded as  special kind of qmm-spaces).
\subsection{Observables and mm-invariants}
Let us recall some of the most important mm-invariants
from \cite{Gro}. Let $Y$ be a compact metric space with
$diam(Y)\leq 1$ and $X\in\chi$. Denote by $Lip_1(X,Y)$ the set of
$1$-Lipschitz functions from $X$ to $Y$.
We can associate to $X$ a compact subset $\cM_Y(X)$ of $\cP(Y)$ by
$$\cM_Y(X)=\{f_*(\mu_X)\mid f\in Lip_1(X,Y)\}\,.$$
One can think about $\cM_Y(X)$ as the information
$Y$ can {\it see} by screening $X$. We will extend
the notion of Lipschitz maps to qmm-spaces and prove that
$\cM_Y(\zeta)$ is well-defined on $\overline{\chi}.$

Let $\nu\in\cP(Y)$ and $0<\kappa<1$.
Then
$diam(\nu,\kappa)$ is the infimal $D$ such that
$D$ contains a subset $Y_0$ such that
$diam(Y_0)\leq D$ and $\nu(Y_0)\geq 1-\kappa$.
The {\it observational diameter} is defined
as
$$Obs Diam_Y(X,\kappa)=\sup_{f\in Lip_1(X,Y)} diam(f_\star(\mu_X),\kappa)\,.$$

\refstepcounter{exa}\label{r5}
\vskip 0.2in
\noindent
{\bf Example} \arabic{exa}
Let $\{K_n\}^\infty_{n=1}$ be the sequence
of complete graphs on $n$ 
vertices. \\ Then $\{\tau(X_{K_n})\}^\infty_{n=1}$ tends to a 
point measure $\alpha$ on $M_\infty$.
Also, $Obs Diam_Y(X_{K_n},\kappa)\to diam(Y,\kappa)$
for any $Y$, $diam(Y)\leq 1/2$ and $0<\kappa<1$.
This follows from the fact that any map $K_n$ to $Y$ is $1$-Lipschitz.
Now, consider the sequence $\{S^n\}^\infty_{n=1}$
of Riemannian spheres of dimension $n$ with diameter $1$ equipped with 
the normalized volume measure.
Then $\{\tau(S^n)\}^\infty_{n=1}$ tends to $\alpha$ as well.
However, $Obs Diam_Y(S^n,\kappa)\to 0$ for any $Y$ and $\kappa$ (by the
L\'evy Concentration Phenomenon). Thus the samplings, in general, do not
capture the observational diameter. However, we have the following
proposition.
\begin{theorem} \label{harmadiktetel} The function $ObsDiam_Y(.,\kappa)$
can be extended to $\overline{\chi}$ in an essentially upper semi-continuous
way. That is, 
for any $Y$ and $0<\kappa'<\kappa $, if $\zeta_m\to\zeta$ weakly then
$$\limsup_{n\to\infty} ObsDiam_Y(\zeta_n,\kappa)\leq ObsDiam_Y(\zeta,\kappa')\,.$$
\end{theorem}
Let $\kappa_1,\kappa_2,\dots,\kappa_N$ be positive numbers such
that $\sum_{i=1}^N \kappa_i<1$.
The separation distance $Sep(X,\kappa_1,\kappa_2,\dots,\kappa_N)$
is the supremal $\delta$ such that there exist Borel sets
$X_i\subset X$, $\mu(X_i)\geq \kappa_i$ such that
$dist_X(X_i,X_j)\geq\delta$.
Obviously, 
$$\lim_{n\to\infty} Sep(X_{K_n},\kappa_1,\kappa_2,\dots,\kappa_N)=1/2,\quad
\,Sep(X_{K_n},\kappa_1,\kappa_2,\dots,\kappa_N)=0\,.$$
\begin{theorem}
\label{negyediktetel}
For any $\kappa_1,\kappa_2,\dots,\kappa_N$,
the function $Sep(.,\kappa_1,\kappa_2,\dots,\kappa_N)$ extends
to $\overline{\chi}$ as an upper semi-continuous function.
\end{theorem}

\section{On the Radon-Nikodym-Dunford-Pettis Theorem}
\label{radon}
Let $(X,\mu,\cA)$ be a probability measure space with a $\sigma$-algebra
$\cA$. Let $f:X\to\bR$ be a bounded measurable function and
$\cB\subset \cA$ be a sub-$\sigma$-algebra.  According to the
Radon-Nikodym Theorem there exists a unique measurable function
$E(f\mid\cB)\in L^\infty(M,\mu,\cB)$ such that
for any $B\in\cB$
$$\int_B E(f\mid \cB)\,d\mu=\int_B f d\mu\,.$$
For the next paragraph our reference is [Chapter 3.]\cite{NM}. 
In our paper we use sometimes Banach valued measurable functions.
In this category, there are several notion of measurability and integral.
What we need is the Gelfand-Dunford integral of weak-$*$-measurable 
functions. So, let $L$ be a Banach space, and $L^*$ be its dual.
An essentially bounded function $f:(X,\mu,\cA)\to L^*$ is 
called weak-$*$-measurable, if for any $v\in L$, the function
$x\to\langle f(x),v\rangle$ is $\cA$-measurable. Note that this is
equivalent to say that as a map $f$ is measurable with respect to
the weak-$*$-topology of $L^*$. We denote the space of these functions
by $l^\infty_{w^*}(X,\mu,\cA,L^*)$.
The Gelfand-Dunford integral $\int f d\mu$ 
of such a function
is the unique element of $L^*$ such that
$$(\int f d\mu)(v)=\int \langle f,v \rangle d\mu$$ holds for any $v\in L$.
Then we still have a Radon-Nikodym type theorem, based on the theorem of 
Dunford and Pettis. The following result was explained to us by Nicolas Monod.
\begin{proposition}[Radon-Nikodym-Dunford-Pettis Theorem]
Let $(X,\mu,\cA)$ and $f$ be as above. Then there exists
an essentially unique function $E(f\mid\cB)$ which is weak-$*$-measurable
with respect to the $\sigma$-algebra $\cB$, such that for all $v\in L$ and 
$B\in\cB$
$$\int_B \langle E(f\mid\cB)(x),v\rangle d\mu = 
\int_B \langle f(x),v\rangle \, d\mu.$$
\end{proposition}
\proof 
By the original Radon-Nikodym Theorem for all $v\in L$ there exists a bounded
$\cB$-measurable function $f_v$ such that for any $B\in\cB$
$$\int_B \langle E(f\mid\cB)(x),v\rangle = 
\int_B f_v d\mu\,.$$
Observe that the map $v\to f_v$ is a continuous linear operator, that
is an element of the space $Hom(L, l^\infty(X,\mu,\cB))$.
Note that $Hom(L,l^\infty(X,\mu,\cB))\sim Hom(l^1(X,\mu,\cB,L^*)$.  On the
other hand, by Proposition 2.3.1 \cite{NM}
$$l^\infty_{w^*}(X,\mu,\cB,L^*)\sim Hom(L,l^\infty(X,\mu,\cB))$$
Hence the function $v\to f_v$ is represented by a weak $*$-measurable
function. \qed

\vskip 0.2in
\noindent
Now let $f:(X,\cA,\mu)\to [-k, k]$ be a
measurable function. Notice that $f$ can be viewed as a function
$\overline{f}:(X,\cA,\mu)\to C[-k,k]^*$, where
$C[-k,k]^*$ is the dual space of the Banach space $C[-k,k]$.
Here $\overline{f}(x)=\delta_{f(x)}$, the point measure concentrated in 
$f(x)$. Then one can consider both
$E(\overline{f}\,\mid\cB)$ and $E(f\,\mid \cB)$. In case of
$f=id:([0,1],\lambda)\to [0,1]$,
$E(f\,\mid\{0,1\})=1/2$ and
$E(\overline{f}\,\mid\{0,1\})$ is the Lebesgue measure.
Now let $Y$ be a compact metric space and 
$f:(X,\cA,\mu)\to Y$ be
a measurable map. Then $E(\overline{f}\,\mid\cB)$ is
a well-defined $\cB$-measurable $C(Y)^*$-valued
function on $X$. On the other hand, in general $E(f\,\mid\cB)$ does
not have a meaning.
If $G:Y\to \bR$ is a continuous function
then $\langle E(\overline{f}\,\mid\cB), G \rangle$ is a measurable
function on $X$ and it is the Radon-Nikodym derivative of 
$G\circ f$. This also shows that
$ E(f\,\mid\cB)$ is a probability measure valued function.

\vskip 0.2in
\noindent
Finally, let us recall the notion of an ultralimit.
Let $Y$ be a compact metric space and $\omega$ be
a nonprincipal ultrafilter. Let $\{y_n\}^\infty_{n=1}\subset Y$
be a sequence of points. Then the ultralimit $\lim_\omega y_n$ is the
unique element $y\in Y$ such that
for any $\e>0$
$$\{n\,\mid d_Y(y_n,y)<\epsilon\}\in\omega $$
Note however that we can define an ultralimit
$\hlimo$ that is valued in $\cP(Y)$, the space of probability metric
spaces on $Y$. Consider the natural
embedding $i:Y\to \cP(Y)$. Then we
have a sequence $\{i(y_n)\}^\infty_{n=1}\subset \cP(Y)$
and the ultralimit $\hlimo\,y_n=\lim_\omega i(y_n)\in\cP(Y)\,.$
Clearly, if $G\in C(Y)$ then
$$\langle \hlimo\, y_n, G \rangle=G(\lim_\omega y_n)\,.$$
We will use the following lemma later.
\begin{lemma} \label{bound}
Let $\cB\subset\cA$ two $\sigma$-algebras on a set $X$ and
let $f:(X,\cA,\mu)\to [0,1]$ be an $\cA$-measurable function.
We denote by $\downarrow E(\overline{f}\mid\cB)(x)$ resp.
by $\uparrow E(\overline{f}\mid\cB)(x)$ the infimum resp. supremum
of the support of the measure $E(\overline{f}\mid\cB)(x)$.
Then if $f\geq\e$ $\mu$-almost everywhere on a set $B\in\cB$, then
$\downarrow E(\overline{f}\mid\cB)(x)\geq \e$ $\mu$-almost everywhere
as well. Similarly, if $\uparrow E(\overline{f}\mid\cB)(x)\leq \e$ 
$\mu$-almost everywhere, 
then $\uparrow E(\overline{f}\mid\cB)(x)\geq \e$ $\mu$-almost everywhere
as well.\end{lemma}
\proof Let $0<\delta<\e$ and $g\in C[0,1]$ such that
$g(t)=1$ if $t\leq \frac{\delta+\e}{2}$, $g(t)=0$ if $t\geq \e$.
Let $C\in\cB$ the subset of $B$ on which
$\downarrow E(\overline{f}\mid\cB)\leq \delta$. Then if $\mu(C)>0$,
$$0<\int_C\langle E(\overline{f}\mid \cB), g\rangle d\mu=
\int_C f\circ g\, d\mu=0\,,$$
leading to a contradiction.

\section{Ultralimits of measured metric spaces}
\label{ultrasection}
Let $\{X_n\}^\infty_{n=1}$ be a sequence of Polish spaces and let
$\omega$ be a nonprincipal ultrafilter. The ultralimit of Polish
spaces was defined in 3.22 \cite{Gro} the following way.
We say that $\{x_i\}^\infty_{i=1}, \{x'_i\}^\infty_{i=1}\subset 
\prod^\infty_{i=1} X_i$ are equivalent if $\lim_\omega d_{X_i}(x_i,x'_i)=0$\,,
where $\lim_{\omega}$ is the associated ultralimit. The elements of
the ultralimit space $\bX$ are the equivalence classes $[\{x_i\}]_{i=1}^\infty$.
The set $\bX$ is a metric space
$$\bd_{\bX}([\bux],[\bux'])=\lim_{\omega} d_{X_i}(x_i,x_i')\,.$$
Then $\bX$ is a complete metric space, but usually it is not separable. 
In order to define a measure on $\bX$ we need some preparation.
 Consider again the spaces
$\{X_i\}^\infty_{i=1}\subset \chi$ and their set-theoretical ultraproduct. 
From now on, we will use the phrase set-theoretical ultralimit, since 
from our point of view, it is very much like the classical ultralimit.
Now, two sequences $\{x_i\}^\infty_{i=1}$ and $\{x'_i\}^\infty_{i=1}$ are equivalent
if $$\{i\,\mid\, x_i=x'_i\}\in \omega\,.$$ The set-theoretical ultralimit is
denoted by $\X$. Then we can define a pseudo-metric on $\X$, by
$$\bd_{\X}({[\bx]},{[\bx']})=\lim_{\omega} d_{X_i}(x_i,x_i')\,.$$
Here, $[\bx]$ denotes the new equivalence class.
Then we have a natural map $\pi:\X\to\bX$. In nonstandard analysis, the inverse
images of an element $[\bux]$ in $\bX$ are called {\it monads} $M([\bux])$. 
The elements of $M([\bux])$ are 
infinitesimally close to each other. Now we define the ultralimit of subsets
the following way:
$[\{a_i\}]\in \A$ if and only if
$$\{i\,\mid\, a_i\in A_i\}\in \omega\,.$$
The ultralimit sets $\A$ form a Boolean algebra $\cP$.
We have a finitely additive measure $\mu_{\X}$
on $\cP$ and this finitely additive measure can be extended to 
a $\sigma$-algebra containing \\ $\cP$ ( see \cite{ESZ} and
\cite{Kechrispaper}) the
following way.
Let $\N\subset\X$ be a {\it nullset} if for any $\e>0$ there exists
an element $\A\in \cP$ such that $\N\subset \A$ and $\mu_{\X}(A)<\epsilon$.
A set $\M\subset\X$ is {\it measurable} if there exists $\P\in\cP$ such
that $\P\triangle \M$ is a nullset. The set of measurable sets $\cM$
is a $\sigma$-algebra with a probability measure, 
where we define $\mu_X(\M)=\mu_X(\P)$. We call a set $\bA\in
\bX$ admissible if $\pi^{-1}(\bA)\in\cM$. The admissible sets form a 
$\sigma$-algebra $\underline{\cM}$. The measure $\mu_{\bX}$ is defined 
on $\bA\in\underline{\cM}$ by $\mu_{\bX}(\bA)=\mu_{\X}(\pi^{-1}(\bA))$.
We do not claim that $\underline{\cM}$ always contains all the Borel-sets, in
the light of the following example. 
\refstepcounter{exa}
\vskip 0.1in
\noindent
{\bf Example} \arabic{exa}. \label{r3}
Let $X_i=K_i$, where $|K_i|=i, d(x,y)=1$ if $x\neq y\in K_i$.
Then $\bX$ is an uncountable discrete set therefore all of its subsets
are Borel-sets. 

\vskip 0.2in
\noindent
However, we prove that balls are always in $\underline{\cM}$, so 
$\underline{\cM}$ contains
the Borel sets, if for some reason $\bX$ is a separable metric space.
\begin{lemma} \label{balls}
If $\bux\in\bX$, $\e\geq 0$, then $B_\e(\bux)\in\underline{\cM}$. Here
$B_\e(\bux)=\{\buy\in\bX\,\mid \bd_{\bX}(\bux,\buy)\leq\e\}$.
\end{lemma}
\proof
We need to prove that $\pi^{-1}(B_\e(\bux))\in\cM$.
Let $\bux=[\{x_i\}]$. Consider the sets
$$\limo B_{\e+\frac{1}{k}}(x_i)=\BB_k\in\cM\,.$$ If
$\by\in \pi^{-1}(B_\e(\bux))$, then $\by\in\cap^\infty_{k=1}\BB_k$. On the other
hand, if $\by\in\cap^\infty_{k=1}\BB_k$ then $\bd_{\X}(\bx,\by)\leq \e$. 
Therefore,
$\pi^{-1}(B_\e(\bux))=\cap^\infty_{k=1}\BB_k\in\cM$.\,\qed

\vskip 0.2in
\noindent
In the course of this paper, we use
bold letters for objects in the set-theoretic ultralimit and
underlined bold letters for the objects in the metric ultralimit.
If $\bA$ is the ultralimit of $\{A_n\}^\infty_{n=1}$, then
we use the notation $\limo A_n=\bA$.

\refstepcounter{exa}
\vskip 0.1in
\noindent
{\bf Example} \arabic{exa}. \label{r4}
For an mm-space $X$, the distance function is Borel on $X\times X$.
This is not always the case for the ultralimit spaces. It is
possible that $d$ is not even $\cM\times \cM$ measurable on $\X\times \X$.
Note that if one considers the set-theoretic ultralimit of the
spaces $X_i\times X_i$, then it is the same space as $\X\times \X$, however
its algebra of measurable functions $\cM_2$ can be much bigger 
than $\cM\times\cM$. This 
phenomenon can be observed if $X_i=X_{G_i}$, where $G_i$ is
a random graph (each edge is chosen with probability $1/2$). Then
with probability $1$, the distance function on the ultralimit
will not be $\cM\times\cM$-measurable (see \cite{ESZ}). However, 
as we shall see soon, the 
distance function on $\X\times\X$ is always $\cM_2$-measurable.

\vskip 0.2in
\noindent
Let $(X,\mu)$ be an mm-space. Recall that the support of $\mu$
is defined the following way. The point $p\in X$ is {\it not}
in the support of $\mu$ if $\mu(B_\e(p))=0$ for some $\e>0$.
Clearly, the support is a closed set with $\mu(Supp(\mu))=1$. Note
however that for some ultralimit spaces such as in Example \ref{r3},
the support of the measure can be empty.

\section{Analysis on the ultralimit}
In this section we fix a sequence $\{X_i\}^\infty_{i=1}\subset\chi$.
As in the previous section $\X$ denotes their set-theoretic ultralimit
and $(\cM,\mu_X)$ stands for the algebra of measurable sets in $\X$ with
the ultralimit measure.
The results in this section are known in the finite graph setting (see 
Section 5 of \cite{ESZ}).
Let $\{f_i:X_i\to [a,b]\}^\infty_{i=1}$ be measurable functions. Their ultralimits
are defined by
$$\ff([x])=\lim_\omega f_i(x_i)\,.$$
\begin{proposition}\label{p41}
The ultralimit function $\ff$ is $\cM$-measurable and
$$\int_{\X}\ff d\mu_{\X}=\limo \int_{X_i} f_i d\mu_{X_i}\,.$$
Conversely, if $\bg:\X\to [a,b]$ is a $\cM$-measurable
function, then there exists a sequence of functions
$\{f_i:X_i\to [a,b]\}^\infty_{i=1}$ such that
their ultralimit $\mu_{\X}$-almost everwhere equals to $\bg$.
\end{proposition}
\proof
In order to prove that $f$ is measurable, it is enough to see that
$$\ff_{[c,d]}=\{\pp\in\X\,\mid c\leq \ff(\pp)\leq d\}\in \cM\,,$$
for any $[c,d]\subset[a,b]\,.$
Let $\P_n=[\{f^i_{[c-\frac{1}{n}, d+\frac{1}{n}]}\}^\infty_{i=1}]\,.$
Clearly, $\ff_{[c,d]}=\cap^\infty_{n=1} \P_n$, thus
$\ff_{[c,d]}\in\cM$. Now fix $k\geq 1$ and
let $h_i:X_i\to\bR$ be a measurable stepfunction such that
$h_i(x)=\frac{j}{2^k}$, when $\frac{j}{2^k}\leq f_i(x)<\frac{j+1}{2^k}$.
Clearly, $|\hh-\ff|\leq \frac{1}{2^k}$ on $\X$. That is
$$\left|\int_{\X} \hh d\mu_{\X}- \int_{\X} 
\ff d\mu_{\X}\right|\leq \frac{1}{2^k}\,.$$
Also, $|\int_{X_i} h_i d\mu_{X_i} -\int_{X_i} f_i d\mu_{X_i}|\leq \frac{1}{2^k}$.
Observe that the ultralimit function 
$\hh$ can be written as $\sum \frac{j}{2^k} \chi_{\C_j}$, where
$\C_j$ is the ultralimit set of $\{C^i_j\}^\infty_{i=1}$, 
$C^i_j=\{x\in X_i\,\mid h_i(x)=\frac{j}{2^k}\}\,.$

\noindent
Therefore $\int_{\X} \hh d\mu_{\X}=\limo \int_{X_i} h_i d\mu_{X_i}\,.$
Consequently, for any $k\geq 1$,
$$\left|\int_{\X} \ff\,d_{\mu_{\X}}-\limo \int_{X_i} 
f_i d\mu_{X_i}\right|\leq\frac{1}{2^{k-1}}\,.$$
Thus 
$$\int_{\X} \ff \,d_{\mu_{\X}}=\limo  f_i d\mu_{X_i}\,.$$
Now let us prove the converse statement.
Let $\bg_k:\X\to [a-1,b+1]$ be the step function
approximation of $\bg$, that is, $\bg_k=\sum \frac{j}{2^k} \chi_{\C_{j,k}}$,
where $\C_{j,k}=\{\bx\in \X\,\mid \frac{j}{2^k}\leq \bg(\bx)< \frac{j+1}{2^k}$.
By modifying $\bg$ on a set of $\mu_{\bX}$-measure zero, we can suppose
that $\C_{j,k}\in\cP$. 
Note that $|\bg(\bx)-\bg_k(\bx)|\leq\frac{1}{2^k}$ on $\X$,
and $|\bg_k(\bx)-\bg_j(\bx)|\leq \frac{1}{2^j}$ if $k\geq j$.
Let $T^i_{j,k}\subset X_i$ be Borel sets such that
$\C_{j,k}=\limo T^i_{j,k}$.
Let $f^i_k=\sum \frac{j}{2^k} \chi_{T^i_{j,k}}\,.$
We define $E^i_k\subset X_i$ in the following
way
$$E^i_k=\{x\in X_i\,\mid |f^i_k(x)-f^i_j(x)|<\frac{1}{2^{j-1}}\,
\mbox{for all $j\leq k$}\}\,.$$
Then $$S_k=\{i\,\mid \mu_{X_i}(E^i_k)>1-\frac{1}{2^k}\}\in\omega\,.$$
Let $h(i)=\min\{i,\sup\{k\mid i\in S_k\}\}$.
Thus for any $k\geq 1$
\begin{equation} \label{cs4}
\{i\,\mid\,h(i)\geq k\}\in\omega\,.
\end{equation}
Let $g'_i=f^i_{h(i)}$. We claim
that $\bg'=\limo g'_i=\bg$ $\mu_{\X}$-almost everywhere.
Let $$V_l=\{\bx\in\X\mid |\bg'(x)-\bg_l(x)|>\frac{1}{2^{l-2}}\}\,.$$
What we need to prove is that for all $l\geq 1$, $\mu_{\X}(V_l)=0\,.$
Since $g_l=\limo f^i_l$, it is enough to show that for
any $k>0$
\begin{equation} \label{cs5}
\{i\,\mid\,\mu_{X_i}(x\,
\mid g'_i(x)-f^i_l(x)|<\frac{1}{2^{l-1}})>1-\frac{1}{2^k}\}\in\omega
\end{equation}
However, (\ref{cs5}) follows from (\ref{cs4}). \qed

\section{Sampling qmm-spaces}
\subsection{q-samplings}\label{qsam}
In this subsection, we show  how one can extend $\tau$ onto qmm's. 
This is described in \cite{LSZDEC} in a slightly different situation
using somewhat different terminology.
Let $M_n$ be the convex compact space of $n\times n$ real matrices
satisfying $0\leq d_{i,j}\leq 1$, $d_{i,j}=d_{j,i}$ and $ d_{i,i}=0\,.$  
Also, let $\hat{M}_n$ be the convex compact space of $\cP[0,1]$-valued
matrices satisfying $d^*_{i,j}=d^*_{j,i}, d^*_{i,i}=\delta_0\,.$
An element
of $\hat{M}_n$ can be viewed as a probability measure on $M_n$.
For a probability measure $\nu$ on $\hat{M}_n$ one can consider its
barycenter $b(\nu)$.
Let $(X,\mu,d^*)$ be a qmm-space. Then the push-forward of
$\mu\times\mu\times\dots\times \mu$ for the natural
 map $\rho_n:X^n\to\hat{M}_n$ is
a probability measure $\nu_n$ on $\hat{M_n}$. Its barycenter 
$b(\nu_n)=\tau_n((X,\mu,d^*))$ is the $n$-sampling measure of $(X,\mu,d^*)$.
Note that $M_\infty$ is the inverse limit of the spaces $M_n$ and $b(\nu_n)$
is the push-forward of $b(\nu)=\tau((X,\mu,d^*))$ constructed in the
Introduction. One can also look at the measures
$b(\nu_n)$ by taking moments.
Let $g=\{g_{ij}:[0,1]\to\bR\}_{1\leq i\leq j\leq n}$ be a system of
continuous functions. They define a continuous function
$q_g=\prod_{1\leq i\leq j\leq n} g_{ij}$ on $M_n$. 
\begin{lemma}
\label{june22}
$$\int_{X^n} \prod_{1\leq i\leq j\leq n} \langle d(x_i,x_j),g_{ij}\rangle
d\mu(x_1) d\mu(x_2)\dots d\mu(x_n)=\int_{M_n} q_g db(\nu_n)\,.$$
\end{lemma}
\proof We have a continuous map $i_g: \hat{M_n}\to\bR$ defined by
$i_g(\underline{s})=\prod_{1\leq i\leq j\leq n}\langle s_{ij},g_{ij}\rangle\,,$
where $\underline{s}=\{s_{ij}\}_{1\leq i\leq j\leq n}\in \hat{M_n}\,.$
By the definition of the push-forward,
$$\int_{X^n} \prod_{1\leq i\leq j\leq n} \langle d(x_i,x_j),g_{ij}\rangle
d\mu(x_1) d\mu(x_2)\dots d\mu(x_n)=\int_{\hat{M_n}} i_g(\underline{s}) 
d\nu_n(\underline{s})\,.$$
On the other hand, by the definition of the barycenter
$$\int_{\hat{M_n}} i_g(\underline{s}) 
d\nu_n(\underline{s})=\int_{\hat{M_n}}
\prod_{1\leq i\leq j\leq n}\langle s_{ij},g_{ij}\rangle d\nu_n(\underline{s}) =
\int_{M^n} (\prod_{1\leq i\leq j\leq n} g_{ij}) d b(\nu_n)\,\quad\qed$$

\vskip 0.2in
\noindent
Following \cite{LSZDEC}, we denote
$\int_{X_r} \prod_{1\leq i\leq j\leq r} \langle d(x_i,x_j), g_{i,j}\rangle
d\mu^r$ by $t(g,x)$.
Observe that the linear combinations of the functions
$q_g$ are dense in $C(M_r)$. Therefore,
$\tau((X_1,\mu_1,d^*_1))=\tau((X_2,\mu_2,d^*_2))$ if and only if
$t(g,X_1)=t(g,X_2)$ for all $r\geq 1$ and system $g$.

\vskip 0.2in
\noindent
Let $\{X_n\}^\infty_{n=1}\subset\chi$ be mm-spaces.
Let $\{X^k_n \}_{n=1}^\infty\subset\chi$ be their $k$-fold product spaces with the
metric $d_{k}(\underline{x},\underline{y})=\max_{1\leq i\leq k} d(x_i,y_i)$.
We consider the ultralimits
$\X, \X^2 \dots$. As we noted before, the $\sigma$-algebras
of measurable sets $\cM_k$ in $\X^k$ are in general much bigger
than the product $\sigma$-algebras $\cM\times\cM\times\dots\times\cM$.
We will denote by $\cM_k^i\subset\cM_k$ the subalgebra
of sets depending on the $i$-th coordinate. That is $\cM_k^i\sim \cM$.
Also, we will denote by $\cM_k^{i,j}$
the subalgebra of sets in $\cM_k$ depending on the $i$-th and the $j$-th
coordinate. That is $\cM_k^{i,j}\sim \cM_2$. The ultralimit
distance function $\bd$ is measurable function in $(\X^2,\mu_{\X}^2,\cM_2)$, 
hence
we can consider the Radon-Nikodym-Dunford-Pettis derivative
$\bd^\star=E(\bd\,\mid\cM\times\cM):\X\times \X\to\cP[0,1]$\,.

\vskip 0.1in
\noindent
A {\it separable realization} of $\bd^*$ is a measurable map
$\Psi:(\X,\mu_{\X},\bd^*)\to ([0,1],\mu_X,d^*)$, where the Borel measure
$\mu_{X}$ is the push-forward  of $\mu_{\X}$, $\bd^*$ is the pull-back
of the Borel function $d^*$. 
Now we show that $([0,1],\mu_X,d_X^*)$ is always
a qmm-space. In fact, we have the following proposition.
\begin{proposition} \label{qmm}
Let $\{X_n\}^\infty_{n=1}\subset \chi$ be mm-spaces and $([0,1],\mu_X,d_X^*)$
be a separable realization of their ultralimits then
$\limo \tau(X_n)=\tau(([0,1],\mu_X,d^*_X))$.
\end{proposition}
\proof
First, let us fix an $r\geq 1$ and a system of continuous functions
as in Subsection \ref{qsam}. By Proposition \ref{p41},
\begin{equation} \label{ESZ6}
\limo\int_{X_n^r} \prod_{1\leq i,j\leq r} g_{ij}(d_{X_n}(x_i,x_j)) d\mu_{X_n}^r=
\int_{\X^r} \prod_{1\leq i,j\leq r} g_{ij} (\bd(\bx_i,\bx_j))d\mu_{\X}^r\,,
\end{equation}
where $\bd$ is the ultralimit of the functions $d_{X_n}$.
\begin{lemma} \label{integration}
$$\int_{\X^r} \prod_{1\leq i,j\leq r} g_{ij} (\bd(\bx_i,\bx_j))d\mu_{\X}^r=
\int_{\X^r} \prod_{1\leq i,j\leq r}  \langle \bd^*(\bx_i,\bx_j),g_{ij}\rangle
d\mu_{\X}^r$$
\end{lemma}
\proof
By definition,
\begin{equation} \label{szerda1}
E(g_{ij}(\bd(\bx_i,\bx_j))\mid\cM_i\times \cM_j)=
\langle \bd^*(\bx_i,\bx_j),g_{ij}\rangle
\end{equation}
By the Integration Rule [Proposition 5.3]\cite{ESZ}, if
$\bh_{ij}:\X^r\to\bR$ are bounded $\cM_{i,j}$-measurable functions then
$$\int_{\X^r} \prod_{1\leq i\leq j\leq r} \bh_{ij} d\mu_{\X}^r=
\int_{\X^r} \prod_{1\leq i\leq j\leq r} E(\bh_{ij}|\cM\times\cM\times\dots\cM)
d\mu_{\X}^r. $$
Hence by (\ref{szerda1}) the lemma follows. \qed

\noindent
Since $\bd^*$ is measurable on $\cL\times\cL$, we immediately 
obtain from the lemma 
that for any separable realization $([0,1],\mu_X, d^*)$:
$$\limo \int_{X_n^r} \prod_{1\leq i\leq j\leq r} 
g_{ij}(d_{X_n}(x_i,x_j)) d\mu_{X_n}^r=
\int_{[0,1]^r} \prod_{1\leq i\leq j\leq r} 
\langle g_{ij}, d^*(x_i,x_j)\rangle d\mu_X^r\,.$$
Note that $\limo \tau(X_n)$ is well-defined since $\cP(M_\infty)$ 
is a compact metric space.
By the definition of the ultralimit, 
$\limo \tau(X_n)=\kappa$ if for any continuous function $g\in C(M_\infty)$
$$\limo \langle \tau(X_n), g\rangle=\langle \kappa, g\rangle\,.$$
Hence our proposition follows. \qed

\vskip0.2in
\noindent
Note that we never actually checked that the
triangle inequality condition holds for the limit qmm-space
$([0,1],\mu_X,d^*)$. However, let $g^{i,j,k}_\e:M_\infty\to\bR $ be a continuous
function, that equals to $1$ if $d_{ij}+d_{jk}-d_{ik}<\e$.
Since $0=\limo \langle \tau(X_n), g^{i,j,k}_\e\rangle=
 \langle \tau([0,1],\mu_X,d^*), g^{i,j,k}_\e\rangle $ 
one can immediately see that the
measure of ``bad'' triangles is zero. That is $([0,1],\mu,d^*)$ is a qmm.
Thus we proved that for any $\kappa\in\overline{\chi}$, there exists
a qmm-space $X$ such that $\tau(X)=\kappa$.
\subsection{Martingales}
Let $(X,\mu,d^*)$ be a qmm-space. The way we obtained $\tau_n(X)$ 
can be described 
by the following sampling process. First pick $n$ points
in $X$ $\mu$-randomly, independently and then for
any $i<j$ choose $d(i,j)$ randomly
according to the probability measure $d^*(x_i,x_j)$.
If $\mu$ is atomless, then by probability one
we obtain a finite mm-space. For measures with atoms, we get
pseudometric spaces. Now let us fix $r\geq 1$ and
a system $g=\{g_{ij}:[0,1]\to\bR\}_{1\leq i\leq j\leq r}$ as in the 
previous section.
Then we have the following martingale $\{B_m,M_m,\cB_m\}^n_{m=0}$.
The $\sigma$-algebras are the standard Borel algebras $\cB_m$ on $M_m$, for
$1\leq m\leq n$. $\cB_0$ is the trivial algebra. The measures are
$\tau_m(X)$, $1\leq m\leq n$. The function $B_n:M_n\to\bR$ is defined
by
$$B_n(Y)=\frac{1}{n^r}\sum_{\phi} t_{\phi}(g,Y)\,,$$
where the summation is taken for all maps $\phi:[r]\to [n]$ and
$$t_\phi(g,Y)=\prod_{1\leq i\leq j\leq r} g_{ij}
(d_Y(i,j))\,.$$
Recall that $t(g,Y)=\prod_{1\leq i\leq j\leq r} g_{ij}(d_Y(p_i,p_j)) d\mu^r(p)$.
Since $Y$ as a measure space is just $[n]$ with the uniform measure,
$t(g,Y)=B_n(Y)\,.$
Then let 
$$B_m=\frac{1}{n^r}\sum_{\phi} E(t_\phi(g)\,\mid M_m)(Y)\,,$$
Note that this is a standard use of martingales in graph limit theory 
see e.g \cite{LSZDEC}.
Let $B_0=\int_{M_n} B_n d\tau_n=t(g,X)$.
Observe that if $\phi$ does not take the value $m$, then
$$E(t_\phi(g)\mid M_m)=E(t_\phi(g)\mid M_{m-1})$$.
Then we have the inequality
$$|B_m-B_{m-1}|\leq \frac{1}{n^r}\sum_{\phi} |E(t_\phi(g)\mid M_m)- 
E(t_\phi(g)\mid M_{m-1})|\,.$$
Otherwise, $|E(t_\phi(g)\mid M_m)-E(t_\phi(g)\mid M_{m-1})|\leq 
\prod_{1\leq i<j\leq r} \|g_{ij}\|=c_g\,.$
Hence,
$|B_m-B_{m-1}|\leq \frac{1}{n}c_g\,.$ Therefore, using the Azuma Inequality
we obtain the following proposition.
\begin{proposition}\label{azuma}
$$Prob(Y\in M_n\mid t(g,Y)-t(g,X)\geq\e)\leq
2 \exp(\frac{-\e^2 n}{2 c_g})\,.$$
\end{proposition}
By the Borel-Cantelli Lemma we immediately obtain the following
corollary.
\begin{proposition}
\label{cantelli}
For any fixed $r\geq 1$ and system $g$, for $\tau(X)$-almost all 
$\zeta\in M_\infty$
$$t(g,X)=\lim_{n\to\infty} t(g,\zeta_n)\,,$$
where $\zeta_n$ is the pseudometric space on the first $n$ coordinate
on $\zeta$.
\end{proposition}

\vskip0.2in
\noindent
So, if $X$ is atomless, then $\zeta_m\in\chi$ and $\tau(X)=
\lim_{n\to\infty} \tau(\zeta_n)$. If $X$ has atoms, then
$\tau(X)=\lim_{n\to\infty} \tau(\zeta_n')$, where
$\zeta'_n$ is the mm-space associated to $\zeta_n$.
Therefore, $\tau(X)\in\overline{\chi}$ for any qmm-space. This
finishes the proof of Theorem \ref{elsotetel}.
\section{The Reconstruction Theorem}
\subsection{The compactification and qmm-spaces}
In this section, we prove Theorem \ref{masodiktetel}.
\begin{proposition} \label{egyikfele}
Let $\chi_Q$ be the metric space of equivalence classes (under the
pseudo-metric $\ubox$) of qmm-spaces. 
Then $\tau:\chi_Q\to\overline{\chi}$ is a continuous bijection.
In particular, if $\ubox(X,Y)=0$ then $\tau(x)=\tau(y)$.
\end{proposition}
\proof
Let $\Psi:([0,1],\lambda)\to X$ be a measure preserving map.
Then clearly, $t(g,X)=t(g,\Psi^{-1}(X)),$ where
$\Psi^{-1}(X)$ is the induced qmm-structure.
Since Lipschitz-functions are dense in $C(Z)$ for any compact 
metric space $Z$, the
proposition follows from the lemma below.
\begin{lemma}\label{lcs1}
For any $r\geq 1$, system of $K$-Lipschitz functions
 $g=\{g_{ij}\}_{1\leq i\leq j\leq r}$ and $\epsilon>0$ there exists
$\delta=\delta_{r,g,\e}>0$ such that if $\square_1(f_1,f_2)<\delta$ for
two weak $*$-measurable $\cP[0,1]$-valued function on $[0,1]^2$, then
$|t(g,([0,1],\lambda,f_1))-t(g,([0,1],\lambda,f_2))|<\epsilon\,.$
\end{lemma}
\proof
$$|t(g,([0,1],\lambda,f_1))-t(g,([0,1],\lambda,f_2))|=$$ $$=
\left| \int_{[0,1]^r} (\prod_{1\leq i\leq j\leq r}
\langle f_1(x_i,x_j), g_{ij},\rangle-
\prod_{1\leq i\leq j\leq r}\langle f_2(x_i,x_j), g_{ij}\rangle)d\lambda^r 
\right|\leq$$
$$\leq \left|\int_{\underline{x}\in[0,1]^r, (x_i,x_j)\in N(f_1, f_2,\e)}
(\prod_{1\leq i\leq j\leq r}\langle f_1(x_i,x_j), g_{ij}\rangle-
\prod_{1\leq i\leq j\leq r}\langle f_2(x_i,x_j), g_{ij}\rangle)d\lambda^r\right| +$$
$$+\left|\int_{\underline{x}\in[0,1]^r, (x_i,x_j)\notin N(f_1, f_2,\e)}
(\prod_{1\leq i\leq j\leq r}\langle f_1(x_i,x_j)g_{ij}\rangle-
\prod_{1\leq i\leq j\leq r}\langle f_2(x_i,x_j)g_{ij}\rangle)d\lambda^r\right|\,,$$
where
$N(f_1, f_2,\e)=\{(x,y)\in X\times X\,\mid d_{ext}(f_1(x,y),f_2(x,y)\leq\e\}\,.$
For the second term, we have the upper bound
$2 c_g\lambda^r(\underline{x}\in[0,1]^r, (x_i,x_j)\in N(f_1, f_2,\e))\leq 
2c_g{r+1 \choose 2} \lambda^2(N(f_1, f_2,\e))\,.$
To estimate the first term of the right hand side of inequality above,
observe that by the definition of the extended metric $d_{ext}$
$$|\langle f_1(x_i,x_j)-f_2(x_i,x_j), g_{ij} \rangle
\leq K d_{ext}(f_1(x_i,x_j)-f_2(x_i,x_j))\,.$$
Note that for positive numbers $c_i,d_i,1\leq i\leq n$ and $T$:
$$|\prod^n_{i=1}(c_i+T)-\prod^n_{i=1} c_i)|\leq (2T)^n\sup_{1\leq i \leq n}
|c_i|^n\,.$$

\noindent
This gives us the upper bound
$$\left|\int_{\underline{x}\in[0,1]^r, (x_i,x_j)\in N(f_1, f_2,\e)}
(\prod_{1\leq i\leq j\leq r}\langle f_1(x_i,x_j), g_{ij}\rangle-
\prod_{1\leq i\leq j\leq r}\langle f_2(x_i,x_j), g_{ij}\rangle)
d\lambda^r\right|\leq$$ $$\leq
(2K\e\sup |g_{ij}|)^{r+1 \choose 2} \,.$$
This immediately shows that if $d(f_1,f_2)$ is small
enough then
$$|t(g,([0,1],\lambda,f_1))-t(g,([0,1],\lambda,f_2))|<\epsilon\,.$$\qed

\subsection{Random maps}
Let $(X,\mu,d^*)$ be a qmm-space. We can suppose that $X=[0,1]$.
Let us pick a sequence $\{x_n\}^\infty_{n=1}$ of independent $\mu$-random points.
For each pair $(i,j)_{i<j}$ we pick
a real number $d(i,j)$ independently according to the probability
measure $d^*(x_i,x_j)$. That is we pick a $\tau(X,\mu,d^*)$-random
element $\underline{x}$ of $M_\infty$.
Let $X_n$ be the restriction of $\underline{x}$ on $[n]$. Then we have
a natural map $\pi_n:X_n\to [0,1]$ defined by $\pi_n(i)=x_i$. We denote by
$\pi$ the ultralimits of the maps $\pi_n$.
\begin{theorem} \label{random}
The map $\pi:\X\to [0,1]$ is a separable realization with probability $1$.
That is $\pi_\star(\mu_\X)=\mu$ and $(\pi\times\pi)^{-1}(d^*)={\bd}_\X^*$ almost
everywhere.
\end{theorem}
\proof
We have two kind of randomness in our construction. First, the random choice
of $\{x_n\}$, then the choice of $d(i,j)$ according to the law 
$d^*(x_i,x_j)$. The following lemma is about the second kind of randomness.
\begin{lemma}
\label{limeszlemma}
For any choice of $\{x_i\}^\infty_{i=1}$ and
$k>0$, with probability one
$$\lim_{n\to\infty}
\frac{\sum_{i\in A_n,j\in B_n} d(i,j)^k}{|A_n||B_n|}=
\frac{\sum_{i\in A_n,j\in B_n}\langle d^\star(x_i,x_j),t^k\rangle}
{|A_n||B_n|}\,. $$
for all sequences $\{A_n, B_n\subset [n]\}^\infty_{n=1}$,
where $|A_n|,|B_n|\geq \e n$ for some $\e>0$. Here $t^k$ denotes the
$k$-th power of the identity function on $[0,1]$.
\end{lemma}
\proof
First, let us recall the Chernoff inequality.
If $X_1,X_2,\dots X_m$ are independent random variables, taking
values in $[0,1]$ and $\delta>0$ then
\begin{equation}\label{chernoff}
Prob\left(\left|\frac{\sum_{i=1}^m X_i}{m}-\frac{\sum_{i=1}^m E(X_i)}{m}\right|>
\delta\right)
\leq 2 \exp\left(  \frac{-\delta^2m}{2}\right)\,.
\end{equation}
Let us apply (\ref{chernoff}) for a fixed pair $A_n, B_n\subseteq [n],
|A_n|,|B_n|\geq \e n$, and the random variables 
$\{d(i,j)^k\}_{i\in A_n,j\in B_n}$.
We get that
$$Prob\left(\left|\frac{\sum_{i\in A_n,j\in B_n} d(i,j)^k}{|A_n||B_n|}-
\frac{\sum_{i\in A_n,j\in B_n}\langle d^\star(x_i,x_j),t^k\rangle}
{|A_n||B_n|}\right|>\delta\right) \leq 2 \exp
\left(  \frac{-\delta^2\e^2n^2}{2}\right)\,. $$
Therefore the probability that 
$$\left|\frac{\sum_{i\in A_n,j\in B_n} d(i,j)^k}{|A_n||B_n|}-
\frac{\sum_{i\in A_n,j\in B_n}\langle d^\star(x_i,x_j),t^k\rangle}
{|A_n||B_n|}\right| $$
is larger than $\delta$ for at least one such pair
is less than $4^n \exp
\left(  \frac{-\delta^2\e^2n^2}{2}\right)\,.$
Hence our lemma follows from the Borel-Cantelli lemma. \qed

\noindent
\vskip 0.2in
Taking ultralimits we immediately obtain the following proposition.
\begin{proposition}
\label{july18}
For any ${\bf A},{\bf B}\in\cM_1$
$$\int_{\bf A}\int_{\bf B} \langle \hat{\bd}^*_{\X}(\bx,\by),t^k\rangle
d\mu_{\X}^2=
\int_{\bf A}\int_{\bf B} \bd_\X(\bx,\by)^k d\mu_{\X}^2\,,$$
where $\hat{\bd}_X^*=(\pi\times\pi)^{-1}(d^\star)\,.$ Consequently,
$(\pi\times\pi)^{-1}(d^\star)=\bd^\star_{\X}$ almost everywhere.
\end{proposition}
\begin{lemma} \label{szinten}
With probability one, $\pi_{\star}(\mu_{\X})=\mu\,.$
\end{lemma}
\proof
Fix a Borel-set $A\subseteq [0,1]$. By the Law of Large Numbers, with
probability one,
$$\lim_{n\to\infty} \frac{|\{i: x_i\in A\}|}{|n|}=\mu(A)\,.$$
That is, with probability one 
$\pi_{\star}(\mu_{\X})$ and $\mu$ coincide on all dyadic intervals. Therefore,
by Caratheodory''s Theorem the two measures are equal. \qed

\vskip 0.2in
\noindent
Now Theorem \ref{random} follows from Proposition \ref{july18} and Lemma
\ref{szinten}. \qed

\subsection{The Maharam Lemma}
Let $([0,1],\cA,\mu)\subset ([0,1],\cB,\mu)$ be two separable $\sigma$-algebras.
We say that $\cA$ is complemented in $\cB$ if there exists a $\sigma$-algebra
$\cC\in \cB$ such that the generated algebra $(\cA,\cC)$ is dense 
in $\cB$ and the elements 
of $\cC$ are
independent from $\cA$. Note that it means that there
exists a $\sigma$-algebra $([0,1],\cC',\mu')$ and the measure preserving
bijection
$$\Phi:([0,1]\times[0,1],\cA\times\cC',\mu\times\mu')\to ([0,1],\cB,\mu)$$
such that $\Phi^{-1}(\cA)=\cA\times [0,1]$. Maharam (\cite{Mah}, see also
\cite{ESZ}) gave a necessary and sufficient condition for having
a complement. Namely, for any $k>0$ there exists a partition
$S_1\cup S_2\cup\dots\cup S_k=[0,1]$, such that $S_i\in\cB$ and $S_i$ is
independent of $\cA$.
\begin{lemma} \label{maharam}
Let $\pi:(\X,\mu_{\X})\to(X,\mu)$ as above a separable
realization. Then there exists
a Maharam partition for any $k\geq 1$ in $\cM^{\X}_1$.
\end{lemma}
\proof
Let $(x_1,x_2,\dots)$ be a random sequence as above. For each $i\geq 1$,\
pick an element $s(x_i)\subset\{1,2,\dots,k\}$ randomly with uniform
distribution.
Then by the Law of Large Numbers, for each dyadic set $A$ and
$1\leq j \leq k$  $$\lim_{n\to\infty}
\frac{|\{1\leq i\leq n\,\mid\,x_i\in A, s(x_i)=j\}}{n}=\frac{1}{k}
\mu(A)\,,$$
with probability $1$.
Hence the ultralimits of $S^n_j=\{s^{-1}(j)\cap [n]\}$,
$\{\bS_1,\bS_2,\dots,\bS_k\}$ form a Maharam partition. \qed
\subsection{The proof of Theorem \ref{masodiktetel}}
Let $([0,1],\mu_X,d^*_X), ([0,1],\mu_Y,d^*_Y)$ be two
qmm-spaces such that $\tau(([0,1],\mu_X,d^*_X))=\tau(([0,1],\mu_Y,d^*_Y))$.
Let us consider a $\tau(([0,1],\mu_X,d^*_X))$-random element of $[0,1]$.
With probability one, both
$\pi_X:(\X,\mu_{\X},\bd^*_{\X})\to ([0,1],\mu_X,d^*_X)$ and
$\pi_Y:(\X,\mu_{\X},\bd^*_{\X})\to ([0,1],\mu_Y,d^*_Y)$ are separable
realizations. By Lemma \ref{maharam}, we have a separable
$\sigma$-algebra in $\cM_1^{\X}$ that contains Maharam-partitions for
any $k>0$ with respect to both $(\pi_X)^{-1}(\cB_{[0,1]})$ and
$(\pi_Y)^{-1}(\cB_{[0,1]})$, where $\cB_{[0,1]}$ is the Borel algebra.
Therefore there exist measure preserving maps
$\rho_X:([0,1],\mu_Z)\to([0,1],\mu_X),$
$\rho_Y:([0,1],\mu_Z)\to([0,1],\mu_Y)$ such that
$(\rho_X\times\rho_X)^{-1}(d^*_X)$ and
$(\rho_Y\times\rho_Y)^{-1}(d^*_Y)$ coincide.
Note that for any Borel probability measure $([0,1],\mu_Z)$ there exists
a measure preserving map $\rho:([0,1],\lambda)\to([0,1],\mu_Z)$. Hence
Theorem \ref{masodiktetel} follows. \qed

\section{Limits vs. Ultralimits}
\begin{proposition}
\label{ubox}
Let $\{X_n\}^\infty_{n=1}\subset \chi$ converge to $X\in\chi$ in the
$\ubox$-metric. Then the support of their metric ultralimit $\bX$
is isometric to $X$. 
Conversely, if for some sequence $\{X_n\}^\infty_{n=1}\subset \chi$
the support of the ultralimit is isometric to $X\in \chi$, then
$X=\limo X_n$. That is for any $\e>0$
$$\{n\,\mid \ubox(X_n,X)<\e\}\in\omega\,.$$
\end{proposition}
\proof
We start with two lemmas.
\begin{lemma}
\label{sep}
Let $(Y,d)$ be a complete metric space
equipped with a probability measure $\mu$, such
that $supp(\mu)=Y$. Then $Y$ is separable.
\end{lemma}
\proof
Let us consider the set $S_{\e/2}$ of $\e/2$-balls
in $Y$. Take a well-ordering $B_1,B_2\dots$ of $S_{\e/2}$ and
for each ordinal $\alpha$, let $f(\alpha)=\mu(\cup_{\beta<\alpha}B_\beta)$.
Then there is a countable ordinal $\alpha$ such that
$f(\alpha)=f(\gamma)$ for any $\gamma\geq \alpha$.
Therefore, we have countable many balls of radius $\e/2$ such
that the measure of their union is greater or equal than
the measure of the union of any countable subset of $S_{\e/2}$. However,
it means that the $\e$-balls around the centers of these balls
cover the whole set $Y$. \qed

\vskip0.2in
\noindent
Now the first part of the proposition easily follows. By (\ref{ESZ6}),
$$\tau(supp(\mu_{\X}))=\tau(X)\,.$$ Hence by the
reconstruction theorem of Gromov, $X$ and $supp(\mu_{\X})$ are isomorphic.
\vskip0.2in
\noindent
Now let us turn to the converse statement.
Let $(X,\mu,d)$ be an mm-space. with $k$ given points
$x_1,x_2,\dots,x_k$ such that
for some $\e>\delta>0$, $\mu(\cup^k_{i=1} B_{x_i}(\delta))\geq 1-\epsilon$,
where $B_{X_i}(\delta)=\{y\in X\,\mid d(x_i,y)\leq \delta\}.$
For $1\leq j \leq k$, let
$C_j=B_{X_i}(\delta)\backslash \cup^{j-1}_{i=1} B_{x_i}(\delta)\,.$
Then consider the following discrete mm-space $(Y,\mu)$:
$Y=\{y_1,y_2,\dots,y_k,z\}$, 
$\nu(y_i)=\nu(C_i)$, $\nu(z)=1-\sum^k_{i=1}\mu(C_i)$,
$d'(y_i,y_j)=d(x_i,x_j)$, $d'(y_i,z)=1\,.$
\begin{lemma}
\label{h1}
Let $(X,\mu)$,$(Y,\nu)$ as above. Then
$\ubox(X,Y)\leq 3\e$.
\end{lemma}
\proof
Let $\Psi:X\to Y$ be defined the following way.
$\Psi(C_i)=y_i$, $\Psi(X\backslash \cup^k_{i=1} C_i)=z\,.$
Then $\square_1(\Psi^{-1}(d'),d)\leq 3\e\,.$
Indeed, $|\Psi^{-1}(d')-d|\leq 3\e$ on $C_i\times C_j$ and
$\mu(\cup^k_{i=1}C_i)>1-\e\,.$ \qed
\vskip 0.2in
\noindent
The following lemma is trivial.
\begin{lemma}
\label{h2}
For each $m>0$ and $\epsilon>0$, there exists some $\delta>0$
such that if
$(G,\mu)$ and $(H,\nu)$ are discrete mm-spaces on the
same set $\{a_1,a_2,\dots,a_m\}$ such that
for all $1\leq i \leq m$ $|\mu(a_i)-\nu(a_i)|\leq \delta$
and for all $1\leq i,j\leq m$, $|d_G(a_i,a_j)-d_H(a_i,a_j)|\leq\delta$
then $\ubox(G,H)\leq\epsilon$.
\vskip 0.2in
\noindent\end{lemma}
Now let $\{X_n\}^\infty_{n=1}\subset\chi$ be mm-spaces
as in the statement of the proposition. and let $\uZZ$ be the support
of $\bX$. By Lemma \ref{balls} and Lemma \ref{sep}, $\uZZ\in\chi$.
Pick $\e>0$. Let $\bux_1,\bux_2,\dots,\bux_k\in \uZZ$ and $0<\delta<\e$ 
such that
\begin{itemize}
\item $\mu(\cup^k_{i=1} B_{\bux_i}(\delta))>1-\e\,.$
\item For each $1\leq i \leq k$, $\mu(B_{\bux_i}(\delta)\backslash 
B'_{\bux_i}(\delta))=0\,,$
where $B'_{\bux_i}(\delta)=\{\buy\,\mid \bd_{\bX}(\bux,\buy)<\delta\}\,$
\end{itemize}
Now let us consider $\{x^n_1,x^n_2,\dots,x^n_k\}\subset X_n$ such
that $\bux_i=[\{x^n_i\}^\infty_{n=1}]$.
For any fix $q>0$
\begin{equation}
\label{h3}
\{n\,\mid \mu_{X_n}(B_{x^n_i}(\delta))-\mu_{\bX}(B_{x_i}(\delta))|< q\}\in\omega
\end{equation}
Indeed, for any $\delta"<\delta<\delta'$
$$\limo \mu_{X_i}(B_{x^n_i}(\delta))\leq\mu_{\bX}(B_{\bux_i}(\delta'))\,\,\,,
\limo \mu_{X_i}(B_{x^n_i}(\delta))\geq\mu_{\bX}(B_{\bux_i}(\delta"))\,.$$
So (\ref{h3}) follows from our condition on negligibility of
the boundary of the balls.
Similarly, for any $j$,
$$\{n\,\mid \mu_{X_n}(C^n_j)-\mu_{\bX}(\C_j)|< q\}\in\omega\,\,,$$
where $C^n_j=B_{x^n_j}(\delta)\backslash \cup^{j-1}_{i=1} B_{x^n_i}(\delta)$.
By Lemmas \ref{h1} and \ref{h2},
$$\{n\,\mid \ubox(X_n,\uZZ)\leq 4\e)\in\omega\,.$$
Hence the proposition follows. \qed

\vskip 0.2in
\noindent
According the proposition, any class of $\chi$ closed
under taking ultralimits is precompact
in $(\chi,\ubox)$.
\refstepcounter{exa}\label{r6}
\vskip 0.2in
\noindent
{\bf Example} \arabic{exa}.
Let $\{C(n)>0\}^\infty_{n=1}$, $\{D(n)>0\}^\infty_{n=1}$ be
two sequences of integers such that
$\lim_{n\to\infty} D(n)=0\,.$
Then $\chi_{C,D}\subset \chi$ is defined as the set of mm-spaces
for which
$$\mu \{x\in X\,\mid \mu(B_{1/n}(x))\leq C(n)\}\leq D(n)\,.$$
By definition, the support of any ultralimit of a sequence in $\chi_{C,D}$
has full measure, so the ultralimit in also in the class $\chi_{C,D}$
(see \cite{Gro} Section 3$1/2$.14)

\section{Lipschitz functions and their ultralimits}
\subsection{Lipschitz maps and qmm-spaces}
Let $(X,\mu,d_X)$ be a qmm-space and $Y$ be a compact metric space.
We say that a measurable map $f:X\to Y$ is a $1$-Lipschitz, if for
$\mu\times\mu$-almost all pairs $(x_1,x_2)\in X\times X$,
$\downarrow d_X(x_1,x_2)\geq d_Y(f(x_1), f(x_2))\,.$
Note that the metric space structure is still well-defined on $Lip_1(X,Y)$.
Also, we can define the sets $\cM(X,Y)\subset\cP(Y)$.
\begin{lemma} \label{ma1}
If $X_1\sim X_2$ ($\ubox(X_1,X_2)=0)$, then
$Lip_1(X_1,Y)$ and $Lip_1(X_2,Y)$ are isomorphic. Also, the
sets $\cM(X_1,Y)$ and $\cM(X_2,Y)$ are equal.
\end{lemma}
\proof
Let $\Phi:[0,1]\to X_1$ be a measurable map. Then
$Lip_1(\Phi^{-1}(X_1),Y)$ is isomorphic to $Lip_1(X_1,Y)$. Hence our
lemma immediately follows from Theorem \ref{masodiktetel}. \qed

\vskip 0.2in
\noindent
Now let $\{X_n\}^\infty_{n=1}\subset\chi$ and $Y$ be a compact metric space. Let
$f_n:X_n\to Y$ be a sequence of $1$-Lipschitz functions.
\begin{lemma} \label{ma2}
The ultralimit of $\{f_n\}^\infty$, $\ff:\X\to Y$ is $1$-Lipschitz and
$\limo (f_n)_\star (\mu_n)=\ff_* (\mu_{\X})\,.$ \end{lemma}
\proof
Let $\bx_1=\limo x^1_n$, $\bx_2=\limo x^2_n$. Then
$$d_{\X}(\bx_1,\bx_2)=\limo d_{X_n}(x^1_n,x^2_n)\geq 
\limo d_Y(f_n(x^1_n),f_n(x^2_n))=d_Y(\ff(\bx_1),\ff(\bx_2))\,.$$

\vskip 0.2in
\noindent
Now let $g:Y\to\bR$ be a continuous function. Then
$\langle (f_n)_\star (\mu_n), g\rangle=\langle \mu_n,g\circ f_n\rangle$
and $\langle \ff_*(\mu_{\X}),g\rangle =\langle \mu_{\X},g\circ \ff\rangle$ Since
$\limo g\circ f_n=g\circ \ff$, by Proposition \ref{p41},
$$\limo \langle (f_n)_\star (\mu_n),g\rangle=\langle 
\ff_\star (\mu_{\X}),g\rangle\,.$$
That is $\limo (f_n)_* (\mu_n)=\ff_*(\mu_{\X})$\,\,\qed

\vskip 0.2in
\noindent
Let us consider the space $(\X,\mu_{\X},\bd^*)$, where $\X$ and $\bd$
are as above. 
\begin{lemma} \label{ma3}
$\ff$ is still $1$-Lipschitz on $(\X,\mu_{\X},\bd^*)$.
\end{lemma}
\proof
Let us apply Lemma \ref{bound} for the $\cM_2$-measurable function
$\bd_{\ff}(\bx,\by)=\bd(\bx,\by)-|\ff(\bx)-\ff(\by)|\,.$
Since $\bd_{\ff}$ is $\mu_{\X}\times\mu_{\X}$-almost positive, the
function $\downarrow \bd^*-|\ff(\bx)-\ff(\by)|$ is 
$\mu_{\X}\times\mu_{\X}$-almost positive as well. 
Therefore, $\ff:(\X,\mu_{\X},\bd^*)\to Y$
is $1$-Lipschitz. \qed
\begin{lemma} \label{ma4}
Let $(X,\mu,d^*)$ be a separable realization
of $(\X,\mu_{\X},\bd^*)$. Then there exists
$\hat{f}\in Lip_1(X,Y)$ such that $\hat{f}_\star(\mu)=\ff_\star(\mu_{\X})\,.$
\end{lemma}
\proof Let $(X',\mu',(d^*)')$ be a separable realization of
$(\X,\mu_{\X},\bd^*)$ for which $\ff$ is measurable. Note that such 
separable realization exists, since any separable extension of a separable
realization is a separable realization. 
Then for this particular realization the function $\hat{f}$ must exist.
Since $Lip_1(X,Y)$ depends only on the $\ubox$-equivalence class
of $X$, the lemma follows. \qed

\vskip 0.2in
\noindent
We can summarize the previous lemmas in a proposition.
\begin{proposition} \label{ez1}
Let $\{X_n\}^\infty_{n=1}\subset\chi$ converging to a qmm-space $(X,\mu,d^*)$
in sampling. Also let $f_n\in Lip_1(X_n,Y)$ such that $(f_n)_*(\mu_n)\to\nu$.
Then there exists $f\in Lip_1(X,Y)$ such that $(f)_*(\mu)=\nu$.
\end{proposition}
\begin{proposition} \label{ez2}
Let $(X,\mu_X,d^*)$ be a qmm-space, $Z$ be a compact metric space and 
$f\in Lip_1(X,Z)$.
Then there exists $\{X_n\}^\infty_{n=1}\subset\chi$
and $\{f_n\in Lip_1(X_n,Z)\}^\infty_{n=1}$ such 
that $(f_n)_*(\mu_n)\to(f)_*(\mu)\,.$ \end{proposition}
\proof Let $\{X_n\}^\infty_{n=1}$ be the discrete spaces
as in Theorem \ref{random}. We define $f_n$
by $f_n(i)=f(x_i)$.
That is $f_n=f\circ\pi_n$, where $\pi_n:X_n\to X$ is the
natural map. Then
$$\limo (f_n)_*(\mu_n)=\ff_*(\mu_{\X})=(f\circ \pi)_*(\mu_{\X})=f_*(\mu)\,.\quad
\qed\,$$

\subsection{The proof of Theorem \ref{harmadiktetel}}
\begin{proposition}
\label{pma1}
Let $Y$ be a compact metric space and $0<\kappa'<\kappa<1$.
Let $\{\nu_n\}^\infty_{n=1},\nu\in \cP(Y)$ such that
$\nu_n\to \nu$. Then
\begin{itemize}
\item $$\limsup_{n\to\infty} diam(\nu_n,\kappa)\leq diam(\nu,\kappa')\,.$$
\item $$\liminf_{n\to\infty}  diam(\nu_n,\kappa)\geq diam(\nu,\kappa)\,.$$
\end{itemize}
\end{proposition}

\proof First, let $Y_0\subset Y$ be a closed set such that $diam(Y_0)=t$
and $\nu(Y_0)\geq 1-\kappa'$.
Let $Y_\e=\{y\in Y\,\mid d_Y(y,Y_0)\leq \e\}\,.$
Pick a continuous function $g:Y\to [0,1]$ such that $g(Y_{\e/2})=1$, 
$g(Y^c_\e)=0$.
Then
$$\limsup_{n\to\infty} \nu_n(Y_\e)\geq \lim_{n\to\infty}\int_Y g d\nu_n=
\int_Y g d\nu\geq 1-\kappa'\,.$$
Hence, for large enough $n$,
$$diam(\nu_n,\kappa)\leq t+ 2\e\,.$$ Therefore,
$\limsup_{n\to\infty} diam(\nu_n,\kappa)\leq diam(\nu,\kappa')$.
One should note that \\
$\limsup_{n\to\infty} diam(\nu_n,\kappa)\leq diam(\nu,\kappa)$ does
not always hold.

\noindent
Now let $Y_n\subset Y$ be closed subsets such that
$diam(Y_n)\leq  t$ and
$\nu_n(Y_n)\geq 1-\kappa$.
We can suppose, by taking a subsequence, that
$\{Y_n\}^\infty_{n=1}$ converges to a closed set $Y_0\subset Y$
in the Hausdorff topology.
Then, $diam(Y_0)\leq t$. Let $g:Y\to [0,1]$ be a continuous
function such that
$g(Y_{\e/2})=1$ and $g(Y^c_{\e})=0$.
Then
$$\liminf_{n\to\infty} \nu_n(Y_n)\leq \lim_{n\to\infty}
\int_Y g\, d\nu_n = \int_Yg d\nu\leq \nu(Y_\e)\,.$$
Hence $\nu(Y_\e)\geq 1-\kappa$ and
$diam(Y_\e)\leq t +2\e\,.$ That is
$\liminf_{n\to\infty}  diam(\nu_n,\kappa)\geq diam(\nu,\kappa)\,.$ \qed

\begin{lemma} \label{mal4}
Let $(X,\mu,d^*)$ be a qmm-space and
$0<\kappa'<\kappa$. Then
\begin{itemize}
\item
for all $\{X_n\}^\infty_{n=1}\subset \chi$ such that
$\{X_n\}^\infty_{n=1}\stackrel{s}{\to} X$
$$\limsup_{n\to\infty} ObsDiam_Y(X_n,\kappa)\leq Obs Diam_Y(X,\kappa')\,.$$
\item There exists $\{X_n\}^\infty_{n=1}\stackrel{s}{\to} X$
such that
$$\liminf_{n\to\infty} ObsDiam_Y(X_n,\kappa)\geq Obs Diam_Y(X,\kappa)\,.$$
\end{itemize}
\end{lemma}
\proof
Let $\{X_n\}^\infty_{n=1}\stackrel{s}{\to} X$.
Pick $f_n\in Lip_1(X_n,Y)$ such that \\ $diam((f_n)_*(\mu_n),\kappa)\geq
ObsDiam_Y (X_n,\kappa)-\frac{1}{n}\,.$
Let $\{(f_{n_k})_*(\mu_{n_k})\}^\infty_{k=1}\in \cP(Y)$ be
an arbitrary convergence subsequence such that
$$\lim_{k\to\infty} diam((f_{n_k})_* (\mu_{n_k}),\kappa)=
\limsup_{n\to\infty} ObsDiam_Y(X_n,\kappa)\,.$$
Then by Proposition \ref{ez1},
there exists $f\in Lip_1(X,\mu)$ such that
$(f_{n_k})_*(\mu_{n_k})\to f_*(\mu)\,.$
By the previous proposition,
$$\limsup_{n\to\infty} ObsDiam_Y(X_n,\kappa)\leq ObsDiam_Y(X,\kappa')\,.$$

Now let $f\in Lip_1(X,Y)$ such that 
$$diam(f_*(\mu),\kappa)\geq Obs Diam_Y(X,\kappa)-\e\,.$$
By Proposition \ref{ez2}, there exists a sequence
$\{X_n\}^\infty_{n=1}\subset \chi$ and $f_n\in Lip_1(X_n,Y)$
such that
\begin{itemize}
\item $\{X_n\}^\infty_{n=1}\stackrel{s}{\to} X$.
\item $(f_n)_*(\mu_n)\to f_*(\mu)\,.$
\end{itemize}
Therefore
\begin{equation}
\label{EV1}
\liminf_{n\to\infty} ObsDiam_Y(X_n,\kappa)\geq 
\liminf_{n\to\infty} diam((f_n)_*(\mu_n),\kappa)\geq
ObsDiam_Y(X,\kappa)-e\,.
\end{equation}
Since (\ref{EV1}) holds for all $\e>0$, the lemma follows.\qed

\vskip 0,2in
\noindent
Now let us finish the proof of our theorem. Let $0<\kappa'<\kappa"<\kappa$,
and $\{(X_n,\mu_n,d^*_n)\}^\infty_{n=1}$ be a sequence of qmm-spaces
converging to $(X,\mu,d^*)$.
Pick a sequence $\{Z_n\}^\infty_{n=1}\subset\chi$ such that
\begin{itemize}
\item $ObsDiam_Y(Z_n,\kappa'')+\frac{1}{n}\geq ObsDiam_Y(X_n,\kappa)\,.$
\item $\{Z_n\}^\infty_{n=1}\stackrel{s}{\to} X$.
\end{itemize}
By the previous lemma,
$$\limsup_{n\to\infty} ObsDiam_Y(X_n,\kappa)\leq
\limsup_{n\to\infty} ObsDiam_Y(Z_n,\kappa")\leq
ObsDiam_Y(X,\kappa')\quad\qed$$

\section{The proof of Theorem \ref{negyediktetel}}
Let $(X,\mu,d^*)$ be a qmm-space,
$0<\kappa_i<1$, $\sum^m_{i=1} \kappa_i\leq 1$. Then
$Sep(X,\kappa_1,\kappa_2,\dots,\kappa_m)$ is
defined as the supremum of $\delta$'s such that
there exist disjoint measurable subsets 
$\{A_i\}^m_{i=1}$, $\mu(A_i)\geq \kappa_i$, with
$$ess\, inf_{x\in A_i, y\in A_j} \downarrow d^*(x,y)\geq \delta\,.$$
Recall that 
$$ess\,inf_U \downarrow d^*= \inf\{t\,\mid 
\mu(p\in U\mid \downarrow d^*(p)\geq t)>0\}\,.$$
\begin{lemma} \label{ucso1}
If $X_1\sim X_2$ then
$Sep(X_1,\kappa_1,\kappa_2,\dots,\kappa_m)=
Sep(X_2,\kappa_1,\kappa_2,\dots,\kappa_m)\,.$ \end{lemma}
\proof Let $\Psi:[0,1]\to X_1$ be a measure preserving map. Then
it is easy to see that $Sep(\Psi^{-1}(X_1),\kappa_1,\kappa_2,\dots,\kappa_m)=
Sep(X_1,\kappa_1,\kappa_2,\dots,\kappa_m)$. Hence the lemma follows from
Theorem \ref{masodiktetel}. \qed
\begin{lemma} \label{ucso2}
Let $\{X_n\}^\infty_{n=1}\subset \chi$ converging in sampling to a 
qmm-space $X$.
Then
$$\limsup_{n\to\infty} Sep(X_n,\kappa_1,\kappa_2,\dots,\kappa_m)
\leq Sep(X,\kappa_1,\kappa_2,\dots,\kappa_m)\,.$$
\end{lemma}
\proof
By taking a subsequence, we can suppose that
$\lim_{n\to\infty}Sep (X_n,\kappa_1,\kappa_2,\dots,\kappa_m)$ exists.
Let $\{A^1_n,A^2_n,\dots,A^m_n\}$ be disjoint subsets such that
$\mu_n(A^i_n)\geq\kappa_i$ and
$$\inf_{i,j} d_{X_n}(A^i_n,A^j_n)\geq 
Sep (X_n,\kappa_1,\kappa_2,\dots,\kappa_m)-\frac{1}{n}\,.$$
Let $\A^i=\limo A^i_n$. Then $\mu_{\X}(\A^i)\geq \kappa_i$ and
$$\inf_{i,j} \bd(\A^i,\A^j)\geq \lim_{n\to\infty} 
 Sep (X_n,\kappa_1,\kappa_2,\dots,\kappa_m)\,.$$
Let $\cL\subset\cM$ be a separable subalgebra such that
$\A^i\subset\cL$ for all $1\leq i \leq m\,.$ Then by Lemma \ref{bound},
$$ess\,\inf \downarrow d^*(\A^i,\A^j)\geq 
\lim_{n\to\infty} 
 Sep (X_n,\kappa_1,\kappa_2,\dots,\kappa_m)\,.$$ 
Therefore we have a separable realization $Y$ of $X$ such that
$$Sep(Y,\kappa_1,\kappa_2,\dots,\kappa_m)\geq lim_{n\to\infty}
Sep(X_n,\kappa_1,\kappa_2,\dots,\kappa_m)\,.$$
Thus the lemma follows from Lemma \ref{ucso1}.  \qed

\vskip 0.2in
\noindent
\begin{lemma}
Let $(X,\mu,d^*)$ be a qmm-space. 
Then there exist $\{Y_n\}^\infty_{n=1}\subset\chi$
such that
$\{Y_n\}^\infty_{n=1}$ converges to $X$ in sampling and
$$\lim_{n\to\infty} Sep(Y_n,\kappa_1,\kappa_2,\dots,\kappa_m)=
Sep(X,\kappa_1,,\kappa_2,\dots,\kappa_m)\,.$$
\end{lemma}
Let $A_1, A_2,\dots, A_m\subset X$ such that
for any $1\leq i<j\leq m$
$$ess\,inf_{x_1\in A_i, x_2\in A_j} \downarrow d^*
(x_1,x_2)\geq t\,.$$
Choose a random sequence $y_1, y_2,\dots$ such that
$$\lim_{n\to\infty} \frac{|(1\leq j\leq n\,\mid y(j)\in A_i)|}{n}=\mu(A_i)\,,$$
and pick $d(i,j)$ randomly according to the law $d^*(y(i),y(j))$.
Let us normalize the  measure on $Y_n$ in such a way that
$\mu_n(B^n_i)=\mu(A_i)$, where
$B^n_i=\{i\,\mid y_i\in A_i\}\,.$
Then with probability one, $(Y_n,\mu_n)$ converges to $X$ in sampling and
for any $1\leq i< j\leq n$, \\
$\inf_{i,j} d_{Y_n}(B^n_i,B^n_j)\geq t$. \qed

\noindent
Now we finish the proof of Theorem \ref{negyediktetel}.
Let $\{X_n,\mu_n,d^*_n)\}^\infty_{n=1}$ be a sequence of qmm-spaces converging
to $(X,\mu,d^*)$ in sampling.
Let us pick a sequence
$\{Z_n\}^\infty_{n=1}\subset\chi$ such that
\begin{itemize}
\item $Sep(Z_n,\kappa_1,\kappa_2,\dots,\kappa_m)+\frac{1}{n}\geq 
Sep(X_n,\kappa_1,\kappa_2,\dots,\kappa_m)\,.$
\item $\{Z_n\}^\infty_{n=1}\stackrel{s}{\to} X$.
\end{itemize}
By the previous lemma,
$$\limsup_{n\to\infty} Sep(X_n,\kappa_1,\kappa_2,\dots,\kappa_m)\leq
\limsup_{n\to\infty} Sep(Z_n,\kappa_1,\kappa_2,\dots,\kappa_m)\leq
Sep(X,\kappa_1,\kappa_2,\dots,\kappa_m)\quad\qed$$

\vskip 0.2in
\noindent
elek@renyi.hu
\end{document}